\documentclass[12pt]{amsart}
\usepackage[utf8]{inputenc}
\usepackage{setspace}
\usepackage[english]{babel}
\usepackage{mathabx}
\usepackage{framed}
\usepackage{soul}
\usepackage[all]{xy}
\usepackage[margin=1in]{geometry}
\usepackage{arydshln}
\usepackage[colorinlistoftodos]{todonotes}
\usepackage{pb-diagram}
\usepackage[mathscr]{euscript}
\usepackage{multicol}
\usepackage{stmaryrd} 
\usepackage{empheq} 
\usepackage{tikz-cd}
\tikzset{
  symbol/.style={
    draw=none,
    every to/.append style={
      edge node={node [sloped, allow upside down, auto=false]{$#1$}}}
  }
}

\usepackage{amsmath,amsfonts,amssymb,amsthm,mathrsfs,latexsym,mathtools}
\usepackage{commath}
\usepackage[T1]{fontenc}
\usepackage{hyperref}
\usepackage{caption} 
\usepackage{subcaption} 

\usepackage{lipsum}

\DeclareMathAlphabet{\mathbbmsl}{U}{bbm}{m}{sl}

\usetikzlibrary{matrix}
\title{Topological Abel-Jacobi map and mixed Hodge structures}

\newcommand{\C}{\mathbb{C}} 
\newcommand{\Z}{\mathbb{Z}}
\newcommand{\R}{\mathbb{R}}
\newcommand{\Hom}{\textup{Hom}}
\newcommand{\Ext}{\textup{Ext}}

\newcommand{\Gr}{\mathrm{Gr}}
\newcommand{\Hv}{H_{\textup{van}}}
\newcommand{\Hp}{H_{\textup{prim}}}
\newcommand{\prim}{\textup{prim}}

\newcommand{\Res}{\textup{Res}}
\newcommand{\Pic}{\textup{Pic}}
\newcommand{\Div}{\textup{Div}}
\newcommand{\Cone}{\textup{Cone}}

\usepackage{amscd} 

\usepackage{enumitem} 

\newtheorem{theorem}{Theorem}[section]
\newtheorem{definition}[theorem]{Definition}
\newtheorem{proposition}[theorem]{Proposition}
\newtheorem{corollary}[theorem]{Corollary}
\newtheorem{question}[theorem]{Question}

\newtheorem{lemma}[theorem]{Lemma}

\newtheorem{remark}[theorem]{Remark}

\newtheorem{notation}[theorem]{Notation}

\newtheorem{thm}{Theorem}

\tikzset{commutative diagrams/.cd,
mysymbol/.style={start anchor=center,end anchor=center,draw=none}
}
\newcommand\MySymb[2][(-1)]{%
  \arrow[mysymbol]{#2}[description]{#1}}

\author{Yilong Zhang}
\address{Department of Mathematics\\
  Purdue University\\
  West Lafayette, IN, 47907, USA}
\email[Y.~Zhang]{zhan4740@purdue.edu}

\date{March 9, 2023}

\begin{document}
\maketitle

\begin{abstract}
    For a smooth projective variety $X$ of dimension $2n-1$ over complex field, Zhao defined the topological Abel-Jacobi map, which sends vanishing cycles on a smooth hyperplane section $Y$ to the middle dimensional primitive intermediate Jacobian of $X$. It agrees with Griffiths' Abel-Jacobi map on vanishing cycles that are algebraic and varies holomorphically on the locus of Hodge classes as hyperplane section deforms. On the other hand, Schnell proposed an alternative construction using the $\R$-split property of the mixed Hodge structure on $H^{2n-1}(X\setminus Y)$. We show that the two definitions coincide, which answers a question of Schnell. 
\end{abstract}

\section{Introduction}

Let $X$ be a compact Riemann surface, and let $J(X)=H^{0}(X,\Omega_X)^{\vee}/H_1(X,\Z)$ be the Jacobian variety of $X$. The Abel-Jacobi map for $X$ is a group homomorphism
\begin{equation}
  A: \Div^0(X)\to J(X)  \label{Intro_Eqn_AJCurve}
\end{equation}
which associates a degree zero divisor $D=\sum_ip_i-q_i$ to a linear functional sending each holomorphic one-form $\omega$ on $X$ to an integral
\begin{equation*}
   \sum_i\int_{q_i}^{p_i}\omega \mod \Lambda, \label{Intro_Eqn_IntegralCurve} 
\end{equation*}
 where $\Lambda$ is the set of values $\int_{\gamma}\omega$, with $\gamma\in H_1(X,\Z)$, called the periods.
 
Since the Abel-Jacobi map is constant on each rational equivalent class, \eqref{Intro_Eqn_AJCurve} factors through $\Div^0(X)/\textup{PDiv}^0(X)\cong \Pic^0(X)$ and induces an abelian group isomorphism 
\begin{equation}
    \bar{A}:\Pic^0(X)\cong J(X).\label{Intro_Eqn_AJCurveIso}
\end{equation}

The injectivity is known as Abel's theorem, while the surjectivity is due to the Jacobi inversion theorem. Moreover, the group $\Pic^0(X)$ carries a natural complex torus structure $H^1(X,\mathcal{O})/H^1(X,\Z)$, called the Picard torus. Under the Dolbeault isomorphism $H^1(X,\mathcal{O})\cong H^{0,1}(X)$, it turns out that the isomorphism \eqref{Intro_Eqn_AJCurveIso} is compatible with the isomorphism of complex tori \cite[Section 4.5]{Clemens-Scrapbook}

\begin{equation}
    \frac{H^{1}(X,\C)}{F^1H^1(X,\C)+H^1(X,\Z)}\cong \frac{H^{1,0}(X)^{\vee}}{H_1(X,\Z)} \label{Intro_Eqn_ComplexTorusIso}
\end{equation}
induced from the unimodular intersection pairing on $H^1(X,\C)$.

From a different perspective, the Picard torus parameterizes the equivalence classes of extension of mixed Hodge structures \cite{Carlson}
\begin{equation}
     0\to H^1(X)\to E\to \Z(-1)\to 0.\label{Intro_Eqn_subseqCurve}
\end{equation}

In particular, a degree zero divisor $D$ determines a mixed Hodge structure $E_D$ fitting into an exact sequence as above. In fact $E_D$ is a sub-mixed Hodge structure of $H^1(X\setminus |D|)$. 
 \begin{proposition}(cf. Proposition \ref{Prop_CurveDoubleInt})\label{Intro_Prop_Curve}
 Under the isomorphism \eqref{Intro_Eqn_ComplexTorusIso}, the  Abel-Jacobi image $A(D)$ corresponds to the extension class $[E_D]$ of an exact sequence of the form \eqref{Intro_Eqn_subseqCurve}.
\end{proposition}

\subsection{Griffiths' Abel-Jacobi Map}
The same result generalizes to the higher dimensions. For a smooth projective variety $X$ of dimension $n$, let $\mathcal{Z}^r(X)_{\textup{hom}}$ denote the group of algebraic cycles of codimension $r$ on $X$ that are homologous to zero. The Griffiths' Abel-Jacobi map \cite{Griffiths} sends a class $Z\in \mathcal{Z}^r(X)_{\textup{hom}}$ to the intermediate Jacobian 
\begin{equation}
J^{2r-1}(X)=\frac{(F^{n-r+1}H^{2n-2r+1}(X,\C))^{\vee}}{H_{2n-2r+1}(X,\Z)}
\end{equation}
by sending $Z$ to a linear functional 
\begin{equation}\label{Griffiths'AJ} 
   [\omega]\mapsto \int_{\Gamma}\omega, 
\end{equation}
where $\Gamma$ is a smooth $(2n-2r+1)$-chain such that $\partial \Gamma=Z$, and $\omega$ is a closed form representing a class in $F^{n-r+1}H^{2n-2r+1}(X,\C)$.

On the other hand, we can consider the trivial Hodge structure $\Z$ with weight $0$, generated by the class $\alpha=[Z]$ and an exact sequence of mixed Hodge structures
\begin{equation}
    0\to H^{2r-1}(X,\Z)\to E\xrightarrow{r} \Z(-r)\to 0,\label{CarlsonAJ}
\end{equation}
where $E$ is a sub-mixed Hodge structure of $H^{2r-1}(X\setminus |Z|)$, and $|Z|$ is the support of $Z$. The extension class of \eqref{CarlsonAJ} is well-defined in $\Ext^1_{\textup{MHS}}(\Z(-r),H^{2r-1}(X,\Z))\cong J^{2r-1}(X)$. This is called Carlson's Abel-Jacobi map.

\begin{proposition} (cf. \cite[Proposition 3.3, 3.4]{KLM})\label{Intro_Prop_CarGri}
Carlson's Abel-Jacobi map coincides with Griffiths' Abel-Jacobi map, up to a sign.
\end{proposition}
The proofs originated from \cite{Jan} and \cite{EV}, using duality between Deligne cohomolgy and Deligne homology. When $Z$ has codimension one, one also refers to \cite[Vol I, Proposition 12.7]{Voisin1} for an equivalent statement and a geometric proof.

\subsection{Topological Abel-Jacobi Map}

In \cite{Zhao}, Zhao introduced the \textit{topological Abel-Jacobi map}. It generalizes Griffiths' Abel-Jacobi map to topological cycles. 

The idea is the following. Let $X$ be a general quintic threefold. It has 2875 rigid lines. Take a hyperplane section $Y$ containing two disjoint lines $L_1, L_2$. The difference $L_1-L_2$ is an algebraic one-cycle homologous to zero, and the integral $\int_{L_2}^{L_1}$ defines Griffiths' Abel-Jacobi map. However, when $Y$ deforms, the algebraic cycle $L_1-L_2$ is obstructed, but the topological class $[L_1]-[L_2]$ deforms with $Y$. Such a class lives in the vanishing cohomology $\Hv^2(Y,\Z)$ \eqref{vanCoh}. The topological Abel-Jacobi map is designed to send these classes to the intermediate Jacobian of $X$ in a continuous way.

Assume $Y$ is smooth and $\alpha\in \Hv^{2}(Y,\Z)$. Let $\gamma$ be a smooth topological 2-cycle representing the Poincar\'e dual class of $\alpha$. There exists a smooth 3-chain $\Gamma$ in $X$ such that $\partial \Gamma=\gamma$. If $\gamma$ is not algebraic, the Griffiths integral \eqref{Griffiths'AJ} does not necessarily vanish on exact $(3,0)+(2,1)$-forms, so it no longer sends $\alpha$ to the intermediate Jacobian. Instead, Zhao's method is to consider harmonic representative $\omega_h$ of a class $[\omega_h]\in F^2H^3(X,\C)$ and add a correction term, namely to consider the sum 
 $$\int_{\Gamma}\omega_h-\int_{\gamma}d^c\sigma$$
where $\sigma$ is a 1-form on $Y$ such that $dd^c\sigma=\omega_{h|_Y}$. Note that this requires $[\omega]$ to lie in the primitive cohomology.

In general, $X$ is a smooth projective variety of dimension $2n-1$. We choose an embedding $X\hookrightarrow \mathbb P^N$ and let $Y$ be a smooth hyperplane section. The topological Abel-Jacobi map is a group homomorphism 
\begin{equation}
   A: \Hv^{2n-2}(Y)\to J_{\prim}(X),\label{TAJ} 
\end{equation}
where $J_{\prim}(X)$ is the intermediate Jacobian \eqref{intJacobian_C} associated with the primitive cohomology $H^{2n-1}_{\prim}(X)$.  The topological Abel-Jacobi map \eqref{TAJ} satisfies the following property \cite[Proposition 2.1.1]{Zhao}:

\begin{enumerate}
\item[\textbf{(P1)}]  \label{P1} If the class $\alpha\in \Hv^{2n-2}(Y,\Z)$ is represented by an algebraic cycle, then $A_{\alpha}$ agrees with Griffiths' Abel-Jacobi image of $\alpha$.
\end{enumerate}

Let $\mathbb P^{sm}$ be the open subspace of projective space parameterizing smooth hyperplane sections of $X$. There is a $\Z$-local system $\mathcal{H}^{2n-2}_{\textup{van},\Z}$ on $\mathbb P^{sm}$ whose fiber at $t$ is the vanishing cohomology $\Hv^{2n-2}(Y_t,\Z)$, where $Y_t$ is the hyperplane section $X\cap H_t$. The underlying \'etale space $T_{\Z}$ of $\mathcal{H}^{2n-2}_{\textup{van},\Z}$ is naturally an analytic covering space of $\mathbb P^{sm}$. Then the topological Abel-Jacobi map \eqref{TAJ} defines a continuous map 
\begin{equation}\label{Eqn_TAJ_family}
\mathcal{A}:T_{\Z}\to J_{\prim}(X).
\end{equation}
Let $\textup{Hdg}\subseteq T_{\Z}$ be the subspace parameterizing pairs $(\alpha,t)$ such that $\alpha\in \Hv^{2n-2}(Y_t,\Z)$, namely the \textit{locus of the Hodge classes}. According to Cattani, Deligne, and Kaplan \cite{CDK}, $\textup{Hdg}$ is a union of algebraic varieties.

The topological Abel-Jacobi map in families \eqref{Eqn_TAJ_family} satisfies the following property \cite[Proposition 2.2.2 and Corollary 2.2.1]{Zhao}:
\begin{enumerate}
\item[\textbf{(P2)}]  \label{P2} The map $\mathcal{A}$ is real analytic.
\item[\textbf{(P3)}]  \label{P3} The map $\mathcal{A}$ is holomorphic on the locus of Hodge classes $\textup{Hdg}$.
\end{enumerate}

As pointed out by Xiaolei Zhao to the author, the property $\textbf{(P3)}$ is predicted by the Hodge conjecture.

\subsection{Construction via Mixed Hodge Structures}
On the other hand, Schnell \cite{Schnell} suggested an alternative way to define a topological Abel-Jacobi map using the $\R$-split property of the mixed Hodge structure on  $H^{2n-1}(X\setminus Y)$ (cf. Corollary \ref{Rsplit_cor}).

Before stating the construction, we first remark that Carlson's Abel-Jacobi map can be defined for Hodge classes as well. Let $U=X\setminus Y$. Then there is a short exact sequence of mixed Hodge structures \eqref{ext_U}
\begin{equation}
    0 \to H^{2n-1}_0(X)\to H^{2n-1}(U)\xrightarrow{\textup{Res}} \Hv^{2n-2}(Y)(-1)\to 0, \label{Intro_Eqn_ext_U}
\end{equation}
which comes from the long exact sequence of cohomologies of the pair $(X,U)$, and we set $H^{2n-1}_0(X):=\textup{Coker}(i_*:H^{2n-3}(Y)\to H^{2n-1}(X))$. The last term is the Tate twist of $\Hv^{2n-2}(Y)$ with the weighted shifted from $2n-2$ to $2n$.

Let $\alpha\in \Hv^{n-1,n-1}(Y,\Z)$ be a Hodge class and let $E_{\alpha}$ be the preimage of $\mathbb Z\alpha$ under the residue map $\textup{Res}$. Then there is an exact sequence 
\begin{equation}
   0\to H^{2n-1}_0(X,\Z)\to E_{\alpha}\to \Z(-n)\to 0 \label{Eqn_HodgeClassExt}
\end{equation}
of mixed Hodge structures as a subsequence of \eqref{Intro_Eqn_ext_U}. The extension class of \eqref{Eqn_HodgeClassExt} lives in $\textup{Ext}^1_{\textup{MHS}}(\Z(-n),H^{2n-1}_0(X,\Z))\cong J_{\prim}(X)$. By functoriality of mixed Hodge structures, this agrees with Carlson's Abel-Jacobi map when $\alpha$ is algebraic.

For $\alpha$ of mixed type, let's forget the complex structure for a moment. Take a real lifting  $\tilde{\alpha}_{\R}\in H^{2n-1}(U,\R)$ and an integral lift $\tilde{\alpha}_{\Z}\in H^{2n-1}(U,\Z)$ of $\alpha$, then $\tilde{\alpha}_{\R}-\tilde{\alpha}_{\Z}$ defines a class in $H^{2n-1}_0(X,\R)$. If we choose a different integral lift $\tilde{\alpha}_{\Z}'$, then the difference $\tilde{\alpha}_{\Z}-\tilde{\alpha}_{\Z}'\in H^{2n-1}_0(X,\Z)$. This defines an image in the real primitive intermediate Jacobian 
\begin{equation}
    J_0(X,\R):=\frac{H^{2n-1}_0(X,\R)}{H^{2n-1}_0(X,\Z)}.\label{Intro_Eqn_J0}
\end{equation}

Of course, the freedom of the choices of $\tilde{\alpha}_{\R}$ forms an affine space that is identified with $H^{2n-1}_0(X,\R)$. But, the $\R$-split property \eqref{conj_eqn} of the mixed Hodge structure $H^{2n-1}(U)$ provides a canonical choice (cf. Proposition \ref{Prop_Rsection}). More precisely, there is a canonical $\R$-linear section $$s_{\R}:\Hv^{2n-2}(Y,\R)(-1)\to H^{2n-1}(U,\R)$$ such that

\begin{enumerate}[label=\text{(\dag)},ref=(\dag)]
\item \label{dag} \textit{$s_{\R}\otimes \C$ is a morphism of $\C$-Hodge structure.} 
\end{enumerate}

The condition \ref{dag} means that $s_{\R}\otimes \C$ sends the $(p,q)$-summand of $\Hv^{2n-2}(Y,\C)$ isomorphically onto $(p+1,q+1)$-summand of $H^{2n-1}(U,\C)$.  For example, when $\dim(X)=1$ as will be discussed in Section \ref{Section_curve}, $s_{\R}$ sends the divisor $\alpha=\sum_ip_i-q_i$ to the unique class in $I^{1,1}_{\R}$ whose residue is $\alpha$, where $I^{1,1}=F^1H^1(U,\C)\cap \overline{F^1H^1(U,\C)}$.

\begin{definition}\normalfont Let $s_{\Z}:\Hv^{2n-2}(Y,\Z)(-1)\to H^{2n-1}(U,\Z)$ be a $\Z$-linear section of \eqref{Intro_Eqn_ext_U}. Call the morphism
\begin{equation}
   \Hv^{2n-2}(Y,\Z)\to J_0(X,\R),~\alpha\mapsto s_{\R}(\alpha)-s_{\Z}(\alpha)\label{realTAJ} 
\end{equation}
Schnell's topological Abel-Jacobi map.
\end{definition}

From a point of view of Carlson's theory (cf. Theorem \ref{thm_Carlson}), the map \eqref{realTAJ} specifies a canonical representive (modulo automorphism over $\Z$) of the extension class of the residue sequence \eqref{Intro_Eqn_ext_U}. One may generalize Schnell's topological Abel-Jacobi map to an arbitrary extension of $\R$-split mixed Hodge structures (cf. Remark \ref{Remark_GeneralTAJ}). 

Note that the condition \ref{dag} guarantees that if $\alpha$ is a Hodge class, $s_{\R}(\alpha)-s_{\Z}(\alpha)$ coincides with the extension class of \eqref{Eqn_HodgeClassExt}. In particular, the map \eqref{realTAJ} satisfies \textbf{(P1)}. Schnell's topological Abel-Jacobi map satisfies the two other properties as well. For \textbf{(P2)}, one notes that $F^{p}$ varies holomorphically, while $\bar{F}^{q}$ varies anti-holomorphically. Therefore, their intersection varies real analytically, and so does $s_{\R}$. \textbf{(P3)} can be proved using variation of mixed Hodge structure.

Christian Schnell asked the following question \cite{Schnell}:

\begin{question} 
Is Schnell's topological Abel-Jacobi map \eqref{realTAJ} the same as Zhao's topological Abel-Jacobi map \eqref{TAJ}?
\end{question}

\subsection{Main Theorem} To make sense of the question, we need to identify the targets of the two maps. Note that there is an isomorphism
$F^{n}\Hp^{2n-1}(X,\C)\cong \Hp^{2n-1}(X,\R), ~a\mapsto \textup{Re}(a)$ between real vector spaces. Next, by the fact that there is a unimodular pairing $\Hp^{2n-1}(X,\Z)\times H^{2n-1}_0(X,\Z)\to \Z$ (cf. Proposition \ref{MHSduality_prop}), the (complex) intermediate Jacobian $J_{\prim}(X)$ is identified with the real Jacobian $J_0(X,\R)$ as a real torus, so the topological Abel-Jacobi map defined by Zhao \eqref{TAJ} can be viewed as a morphism to the real torus. 

The next result answers this question affirmatively:

\begin{theorem} \label{main theorem}
Schnell's topological Abel-Jacobi map and Zhao's topological Abel-Jacobi map coincide. Moreover, the map associates $[\psi]\in \Hv^{2n-2}(Y,\Z)$ with a linear functional
$$ [\omega]\mapsto -\int_X\omega\wedge \varphi+\int_Yd^c\sigma\wedge \psi,\mod \int_{X}\omega\wedge\lambda_{\Z},$$
where $\lambda_{\Z}\in \Hp^{2n-1}(X,\Z)$, $\omega$ is a closed form of pure type $(p,q)$ representing a class $[\omega]\in F^n\Hp^{2n-1}(X,\C)$ and $\sigma$ is a form on $Y$ such that $\omega_{|Y}=dd^c\sigma$. $\varphi$ is a smooth $(2n-1)$-form on $X$ with a log pole along $Y$ such that $\mathcal{R}es_Y(\varphi)=\psi$.
\end{theorem}

As an additional remark, the topological Abel-Jacobi map can be defined on a hyperplane section $Y_{nd}$ that has an ordinary double point. This is because the cohomology $H^{2n-2}(Y_{nd})$ is still a pure Hodge structure, and the mixed Hodge structure on $H^{2n-1}(X,Y_{nd})$ is still $\R$-split.

To prove the theorem, we first interpret Zhao's topological Abel-Jacobi map as a topological pairing between relative cohomology and relative homology. Next, we show the map $\omega\mapsto (\omega, d^c \sigma)$ defines the canonical $\R$-splitting of the mixed Hodge structure sequence \eqref{ext}
\begin{equation}
    0\to H^{2n-2}_0(Y)\to H^{2n-1}(X,Y)\to \Hp^{2n-1}(X)\to 0. \label{Intro_Eqn_ext}
\end{equation}

Finally, by Poincar\'e duality and Lefschetz duality, we show the sequence \eqref{Intro_Eqn_ext} is dual to \eqref{Intro_Eqn_ext_U} as extensions of $\R$-split mixed Hodge structures, up to a Tate twist. Then reduce the problem to a linear algebra argument.

\textbf{Structure of the paper.} In section \ref{Section_Prelim}, we review definitions of primitive and vanishing cohomologies, $\mathbb R$-split mixed Hodge structures and their extensions. Our main examples are mixed Hodge structures on $H^{2n-1}(X,Y)$ and $H^{2n-1}(U)$. We will study a non-degenerate pairing between them.

In section \ref{Section_curve}, we study the relation between Abel-Jacobi map and extensions of mixed Hodge structures in curve case and provide a proof for Proposition \ref{Intro_Prop_Curve}.

In section \ref{Section_ZhaoTAJ}, we review the definition and properties of Zhao's topological Abel-Jacobi map. We obtain a different formula \eqref{Eqn_ZhaoPairingResidue} for the topological Abel-Jacobi map using differential forms with logarithmic poles along a divisor, which leads to a proof of the main theorem in the special case when $\alpha$ is a Hodge class. In the end, we relate $\partial\bar{\partial}$-lemma to the fact that the mixed Hodge structure $H^{2n-1}(X,Y)$ is $\R$-split.

In section \ref{Section_Proof}, we will prove the Theorem \ref{main theorem}. In Appendix \ref{Section_App}, we aims to find an explicit pairing between $H^{2n-1}(X,Y)$ and $H^{2n-1}(U)$. In Appendix \ref{Sec_AppB}, we will work out the definition of relative cohomolgy and Hodge filtration using hypercohomolgy in the geometric case.

\textbf{Acknowledgement}: I would like to thank my advisor, Herb Clemens, for introducing me to this research project and for his support. I would like to thank Christian Schnell for sharing his idea on the construction of topological Abel-Jacobi map by $\R$ split mixed Hodge structures. I would like to thank Xiaolei Zhao for explaining the details of his thesis to me. I would like to thank anonymous referee for various suggestions to improve the manuscript.

\section{Preliminaries} \label{Section_Prelim}
\subsection{Primitive and Vanishing Cohomologies}
Let $X$ be a smooth projective variety of dimension $2n-1$ and let $Y$ be a smooth hyperplane section with respect to a projective embedding. Let $i:Y\subseteq X$ denote the inclusion. The middle dimensional primitive cohomology $\Hp^{2n-1}(X,\Z)$ on $X$ is defined as 
\begin{equation}
    \Hp^{2n-1}(X,\Z):=\ker (i^*:H^{2n-1}(X,\Z)\to H^{2n-1}(Y,\Z)).\label{primCoh}
\end{equation}

There is an associated primitive intermediate Jacobian
\begin{equation}
 J_{\prim}(X):=(F^n\Hp^{2n-1}(X,\C))^{\vee}/H_{2n-1}(X,\Z)_{\prim}.\label{intJacobian_C}
\end{equation}

Here we take the torsion free part of integral homology whenever necessary. 
The vanishing cohomology $\Hv^{2n-2}(X,\Z)$ is defined as the kernel of Gysin homomorphism
\begin{equation}
    \Hv^{2n-2}(X,\Z):=\ker(i_*:H^{2n-2}(Y,\Z)\to H^{2n}(X,\Z)).\label{vanCoh}
\end{equation}

 Equivalently, it consists of Poincar\'e dual classes of homology classes on $Y$ that go to $0$ in the homology on $X$. See \cite[Vol II, Chapter 2]{Voisin2} for basic properties of primitive and vanishing cohomologies.

Both $\Hp^{2n-1}(X,\Z)$ and $\Hv^{2n-2}(X,\Z)$ are pure Hodge structures of weight $2n-1$ and $2n-2$, respectively.

\subsection{Deligne's $\R$-split MHS} 
Recall that from Deligne's theory \cite{Deligne}, a mixed Hodge structure is a triple $H=(H_{\Z}, W, F)$, where $H_{\Z}$ is a free abelian group, $W_{\bullet}H_{\Z}$ is an increasing weight filtration and $F^{\bullet}H_{\C}$ is a decreasing Hodge filtration such that each graded piece $\Gr^W_iH_{\C}$ is a pure Hodge structure of weight $i$, where $H_{\C}:=H_{\Z}\otimes \C$.

Let $H$ be a mixed Hodge structure. Then one can define the dual mixed Hodge structure $H^{\vee}$ whose weight filtration is
\begin{equation}
    W_kH^{\vee}_{\Z}=\{\varphi\in H^{\vee}_{\Z}|\varphi(W_{-k-1}H_{\Z})=0\}.\label{DualWeight_def}
\end{equation}


The Hodge filtration of the dual space $H^{\vee}$ is defined as 
\begin{equation}
    F^pH^{\vee}_{\C}=\{\varphi\in H^{\vee}_{\C}|\varphi(F^{-p+1}H_{\C})=0\}.\label{DualHodge_def}
\end{equation}

A morphism $f:H\to H'$ between mixed Hodge structures (of weight zero) is a $\Z$-linear map $f_{\Z}: H_{\Z}\to H'_{\Z}$ that preserves weight filtration and its complexification preserves Hodge filtrations.

Every mixed Hodge structure $H$ admits a unique splitting over $\C$ \cite[Theorem 2.13]{CKS}
\begin{equation}\label{eqn_Csplit}
H_{\C}\cong \bigoplus_{p,q}I^{p,q},\ \mbox{such that}\ \ F^pH_{\C}\cong \bigoplus_{k\ge p}I^{k,q},\ W_kH_{\C}\cong \bigoplus_{p+q\le k}I^{p,q},
\end{equation} 
which satisfies the conjugation property 
\begin{equation*}
    \bar{I}^{p,q}\equiv I^{q,p}\mod \bigoplus_{k< p,l< q}I^{k,l}. \label{conj}
\end{equation*}
\begin{definition}\normalfont \cite[p.64]{MHS} 
A mixed Hodge structure $H$ is called \textit{$\R$-split} if the conjugation property of Deligne's decomposition is equality, namely
\begin{equation}
    \bar{I}^{p,q}= I^{q,p},~\forall p,q.\label{conj_eqn}
\end{equation}
\end{definition}
In particular, pure Hodge structures are always $\R$-split.

\begin{proposition}\label{Rsplit_prop}
Suppose the mixed Hodge structure $H$ has weights concentrated on levels $w-1$ and $w$ consecutively, then $H$ is $\R$-split.
\end{proposition}
\begin{proof}
Let $p+q=w-1$ or $w$.  Then for all $k<p$ and $l<q$, $k+l<w-1$ so $I^{k,l}=0$. So the condition \eqref{conj_eqn} trivially holds. 
\end{proof}

\subsection{Extension of Mixed Hodge Structures}

Let $A, B$ be two mixed Hodge structure with $B>A$, namely $W_nA=A$, and $W_nB=0$ for some $n\in \mathbb Z$. An extension of mixed Hodge structures of $B$ by $A$ is an exact sequence
\begin{equation}
   0\to A\xrightarrow{f} E\xrightarrow{g} B\to 0, \label{Eqn_Carlson_ext}
\end{equation}
where $f$ and $g$ are morphisms of mixed Hodge structures, which means they are morphisms of underlying abelian groups and preserve both weight and Hodge filtrations. 

\begin{proposition}\label{Prop_Rsection}
Suppose $A$ and $B$ are $\R$-split mixed Hodge structures, then $E$ is $\R$-split if and only if there is a $\R$-linear section $s_{\R}:B_{\R}\to E_{\R}$ such that $s_{\R}\otimes \C$ preserves Hodge filtrations. Moreover, such a section is unique.
\end{proposition}
\begin{proof}
Denote $E_{\C}\cong\oplus_{p,q}I^{p,q}$, $A_{\C}\cong\oplus_{i,j}K^{i,j}$, and $B_{\C}\cong\oplus_{s,t}J^{s,t}$ the $\C$-splitting as in \eqref{eqn_Csplit}. Then $g$ sends $I^{p,q}$ isomorphically onto $J^{p,q}$ for $p+q> n$. 

Suppose $E$ is $\R$-split. We define a $\C$-linear section $s:B_{\C}\to E_{\C}$ to be the inverse of $g_{|I^{p,q}}$ on each $J^{p,q}$ and extends $\C$-linearly. Then by definition, $s$ preserves Hodge filtrations. It remains to show that $s$ descends to a section over $\R$. Note that the condition \eqref{conj_eqn} implies $B_{\R}=\oplus_{p,q} B_{\R}\cap (J^{p,q}\oplus J^{q,p})$. So it suffices to show that for each $(p,q)$, the image of $B_{\R}\cap (J^{p,q}\oplus J^{q,p})$ live in $E_{\R}$.  For each $b\in B_{\R}\cap (J^{p,q}\oplus J^{q,p})$, it can be expressed as $b=c_{p,q}+\overline{c_{p,q}}$, where $c_{p,q}\in J^{p,q}$. Then by the $\C$-linearity of $s$, we have $s(b)=s(c_{p,q})+\overline{s(c_{p,q})}\in E_{\R}$. 

Conversely, if there is such a section $s_{\R}$, then $s_{\R}(J^{p,q})$ together with $f_{\C}(K^{p,q})$ defines a $\C$-splitting of $E_{\C}$ satisfying \eqref{conj_eqn}. So $E$ is $\R$-split.

For uniqueness, suppose $s_1$ and $s_2$ are two such sections, then $s_1-s_2$ lifts to a map $B_{\R}\to A_{\R}$ and preserves Hodge filtrations after complexification, then $(s_1-s_2)_{\C}$ must send $J^{p,q}$ to $K^{p,q}$, which is zero since $p+q>n$ and $A$ has weight at most $n$. 
\end{proof}
\begin{definition}\normalfont
We call the $\R$-linear section $s_{\R}$ \textit{Deligne's $\R$-linear section}.
\end{definition}

Note in general, $s_{\R}$ is not defined over $\Z$. Otherwise, the extension \eqref{Eqn_Carlson_ext} splits over $\Z$.

Now let's review a theorem of Carlson characterizing equivalent classes $\textup{Ext}^1_{\textup{MHS}}(B,A)$ of extensions of $B$ by $A$ of the form \eqref{Eqn_Carlson_ext} in the category of mixed Hodge structures.

\begin{theorem}(Carlson \cite{Carlson})\label{thm_Carlson}
There is an isomorphism 
\begin{equation*}
    \textup{Ext}^1_{\textup{MHS}}(B,A)\cong J^0\Hom(B,A),
\end{equation*}
where $J^0\Hom(B,A)$ is the complex torus 
$$\frac{\Hom_{\mathbb C}(B,A)}{F^0\Hom_{\mathbb C}(B,A)+\Hom_{\mathbb Z}(B_{\mathbb Z},A_{\mathbb Z})},$$
and $F^0\Hom_{\mathbb C}(B,A)$ consists of $\C$-linear maps $B_{\C}\to A_{\C}$ preserving Hodge filtrations. 

The extension class $[E]$ of \eqref{Eqn_Carlson_ext} is represented by a homomorphism 
\begin{equation}\label{eqn_sF-sZ}
s_F-s_{\Z}:B_{\C}\to A_{\C},
\end{equation}
where $s_F:B_{\C}\to E_{\C}$ is a $\C$-linear section that preserves Hodge filtrations, and $s_{\Z}:B_{\Z}\to E_{\Z}$ is a $\Z$-linear section.
\end{theorem}

\begin{remark}\normalfont \label{Remark_GeneralTAJ}
For an extension of $\R$-split mixed Hodge structures of the form \eqref{Eqn_Carlson_ext} in general, one may extend Schnell's construction \eqref{realTAJ} to define a topological Abel-Jacobi map
$$s_{\R}-s_{\Z}: B_{\Z}\to A_{\R}/A_{\Z}.$$
\end{remark}

\subsection{$\partial\bar{\partial}$-Lemma}
On a compact K\"ahler manifold, the $\partial\bar{\partial}$-lemma \cite[p.149]{GH} states that an exact form of type $(p,q)$ can be expressed as $\partial\bar{\partial}\sigma$ for some form $\sigma$ of type $(p-1,q-1)$.

Let $d^c=i(\bar{\partial}-\partial)$, then both $d=\partial+\bar{\partial}$ and $d^c$ are real operators. Moreover, $dd^c=2i\partial\bar{\partial}$. So up to replacing $\sigma$ by a constant multiple, the exact form can be expressed as $dd^c\sigma$. $\partial\bar{\partial}$-lemma is also known as $dd^c$-lemma, especially for real differential forms.

\begin{proposition}\label{Prop_ddbarXY}
Let $X$ and $Y$ be the same as at the beginning of Section \ref{Section_Prelim}. Let $\omega$ be a closed $(p,q)$-form on $X$ representing a class in the primitive cohomology $\Hp^{2n-1}(X,\C)$, then there is a $(p-1,q-1)$-form $\sigma$ on $Y$ such that 
$$\omega_{|Y}=dd^c\sigma.$$ 
\end{proposition}
\begin{proof}
According to the definition of primitive cohomology \eqref{primCoh}, the restriction $[\omega_{|Y}]$ is zero in the cohomology of $Y$. Therefore $\omega_{|Y}$ is an exact form. So the claim follows from $\partial\bar{\partial}$-Lemma.
\end{proof}

\subsection{Relative Cohomology}

Let $X$ be a smooth projective $2n-1$ fold, and let  $i:Y\to X$ be the inclusion of a smooth hyperplane section. Denote $U:=X\setminus Y$ the complement.

The relative cohomology $H^{2n-1}(X,Y)$ modulo torsion can be defined through de Rham theory: For a smooth manifold $M$, denote $A^{\bullet}_M=\{\cdots\to A_M^k\xrightarrow{d_M} A_M^{k+1}\to \cdots\}$ the complex of $C^{\infty}$ forms on $M$. Then the relative de Rham cohomology $H_{dR}^{*}(X,Y)$ is the cohomology theory associated to the mapping cone complex of $i^*:A^{\bullet}_X\to A^{\bullet}_Y$ (cf. \cite[p.78]{BottTu}). The differential at $k$-th term is 
$$d:A^k_X\oplus A^{k-1}_Y\to A^{k+1}_X\oplus A^k_Y$$
$$d(\omega,\tau)=(d\omega,\omega_{|Y}-d\tau).$$ 
This defines relative cohomology with $\R$ coefficients, and similarly for $\C$ coefficients, if we take differential forms to be complex valued. To define integral structure of $H_{dR}^{*}(X,Y)$, we take the classes whose pairing to singular relative homology with $\Z$ coefficients $H_*(X,Y,\Z)$ being integral valued. 

Alternatively, one can use the hypercohomology of holomorphic de Rham complex to define the relative cohomology. The benefit is that the torsion information is remembered, and the compatible integral structure and the Hodge filtration can be defined through resolution and filtration on complex of sheaves. In Appendix \ref{Sec_AppB}, we will relate the two definitions. 

The relative cohomology $H^{*}(X,Y)$ carries a mixed Hodge structure making the long exact sequence for the pairs $(X,Y)$
\begin{equation}
    \cdots \to H^{2n-2}(X) \xrightarrow{i^*} H^{2n-2}(Y)\to H^{2n-1}(X,Y)\to H^{2n-1}(X)\xrightarrow{i^*} H^{2n-1}(Y)\to \cdots, \label{les_XY}
\end{equation}
an exact of mixed Hodge structures \cite[Theorem 3.22]{MHS}. The Hodge filtration $F^pH^{\bullet}(X,Y,\C)$ is the subspace represented by closed forms in $F^pA^{\bullet}_X\oplus F^pA^{\bullet-1}_Y$ (cf. Proposition \ref{Prop_App_Fp}).

The weight filtration $W_kH^{*}(X,Y,\Z)$ can be read off from the exact sequence \eqref{les_XY} and the fact that both of the kernel and cokernel of $i^*$ are all pure Hodge structures. In fact, by splitting up the long exact sequence \eqref{les_XY}, we obtain a short exact sequence  
\begin{equation}
   0\to H^{2n-2}_0(Y)\to H^{2n-1}(X,Y)\to \Hp^{2n-1}(X)\to 0,\label{ext}
\end{equation}
where $H_0^{2n-2}(Y,\Z):=\textup{coker} (i^*:H^{2n-2}(X,\Z)\to H^{2n-2}(Y,\Z))$ carries pure Hodge structure of weight $2n-2$. So $H^{2n-1}(X,Y)$ has weight in $2n-2$ and $2n-1$.
\subsection{Cohomology of Hyperplane Complement}
The cohomology $H^{2n-1}(U)$ on the complement also admits a canonical mixed Hodge structure due to Deligne. $H^{k}(U)$ is the $k$-th hypercohomology of the log complex $\Omega^{\bullet}_X(\log Y)$. The Hodge filtration $F^pH^{k}(U,\C)$ is the $k$-th hypercohomology of the subcomplex $\Omega^{\ge p}_X(\log Y)$. The log complex admits an acyclic resolution using the double complex $\oplus_{p,q}\mathcal{A}_X^{p,q}(\log Y)$ (cf. \cite[Vol I, p.212]{Voisin1}). Let $A_X^{p,q}(\log Y)$ be the space of global sections of $\mathcal{A}_X^{p,q}(\log Y)$, then $F^pH^{k}(U,\C)$ consists of classes that are represented by closed forms in $\oplus_{l\ge p}A_X^{l,k-l}(\log Y)$. Namely, they are complex valued $C^{\infty}$ $k$-forms with a log pole along $Y$ and have at least $p$ holomorphic differentials.

\begin{remark}\normalfont
All cohomology are over the integral ring, where we omit the coefficients symbol $\Z$. Moreover, we always take the torsion-free part of cohomology groups.
\end{remark}

The weight filtration $W_kH^*(U,\Z)$ is determined by the exact sequence of mixed Hodge structures coming from the long exact sequence of cohomologies of the pair $(X,U)$
\begin{equation}
    \cdots \to H^{2n-3}(Y) \xrightarrow{i_*} H^{2n-1}(X)\to H^{2n-1}(U)\xrightarrow{\Res} H^{2n-2}(Y)(-1)\xrightarrow{i_*} H^{2n}(X)\to \cdots, \label{les_U}
\end{equation}
where we used the Thom isomorphism $H^{k}(X,U)\cong H^{k-2}(Y)(-1)$ which is an isomorphism of mixed Hodge structures. $i_*$ is the Gysin homomorphism.

Similarly, by splitting up \eqref{les_U}, we obtain the short exact sequence
\begin{equation}
  0\to  H^{2n-1}_0(X)\to H^{2n-1}(U)\xrightarrow{\textup{Res}} \Hv^{2n-2}(Y)(-1)\to 0. \label{ext_U}
\end{equation}

\begin{corollary}\label{Rsplit_cor}
Both $H^{2n-1}(X,Y)$ and $H^{2n-1}(U)$ are $\R$-split mixed Hodge structures.
\end{corollary}
\begin{proof}
From the exact sequence \eqref{ext}, the mixed Hodge structure of $H^{2n-1}(X,Y)$ has weight concentrated on levels $2n-2$ and $2n-1$ consecutively. Similarly, by the exact sequence \eqref{ext_U}, $H^{2n-1}(U)$ has weights concentrated on levels $2n-1$ and $2n$ consecutively. Then the claim follows from Proposition \ref{Rsplit_prop}.
\end{proof}

\subsection{Lefschetz Duality}
Note that the sequence is isomophic to homology sequence 
\begin{equation}\label{eqn_lesHomology}
\cdots \to H_{2n-1}(Y) \xrightarrow{} H_{2n-1}(X)\to H_{2n-1}(X,Y)\xrightarrow{} H_{2n-2}(Y)\xrightarrow{} H_{2n-2}(X)\to \cdots
\end{equation}
using the Poincar\'e duality and the Lefschetz duality \cite[Theorem B.28]{MHS}
\begin{equation}\label{Eqn_LD}
H^{2n-1}(U,\Z)\cong H_{2n-2}(X,Y,\Z).
\end{equation} 

The natural pairing between \eqref{les_XY} and \eqref{eqn_lesHomology} induces a pairing between \eqref{les_XY} and \eqref{les_U}, and then induces a pairing between the short exact sequences \eqref{ext} and \eqref{ext_U}:

\begin{figure}[h]
    \centering
    \begin{equation}\label{diagram}
\begin{tikzcd}
0\arrow[r]&H^{2n-2}_0(Y,\Z)\MySymb{dr}\arrow[d,symbol=\times] \arrow[r,"f"] &  H^{2n-1}(X,Y,\Z) \arrow[d,symbol=\times]\arrow[r,"g"]&\Hp^{2n-1}(X,\Z)\arrow[d,symbol=\times]\arrow[r]&0\\
0&\Hv^{2n-2}(Y,\Z)(-1) \arrow[l] \arrow[d]& H^{2n-1}(U,\Z)\arrow[l,"q"'] \arrow[d]&H^{2n-1}_0(X,\Z)\arrow[l,"p"'] \arrow[d]&0\arrow[l]\\
&\Z&\Z&\Z&
\end{tikzcd}
\end{equation}
\end{figure}
The pairings on the first and third columns are induced from the intersection pairing on $H^{2n-2}(Y,\Z)$ and $H^{2n-1}(X,\Z)$ given by wedge product, respectively.

The pairing
\begin{equation}
    H^{2n-1}(X,Y,\Z)\times H^{2n-1}(U,\Z)\to \Z\label{RelativeOpenPair}
\end{equation}
in the middle is non-degenerate (modulo torsion) and is explicit as follows (cf. Proposition \ref{ExplicitPairing_Prop}):
\begin{equation}
    \big( [\omega,\tau], [\varphi] \big)\mapsto \int_X\omega\wedge \varphi-\int_{Y}\tau\wedge \mathcal{R}es_Y(\varphi).\label{RelativeOpenPair_Formula}
\end{equation}

Here $\varphi$ is a $C^{\infty}$ form with a log pole along $Y$, namely, it can be locally expressed as $\eta\wedge\frac{dz}{z}+\eta'$ for some $C^{\infty}$ forms $\eta,\eta'$. So the integral $\int_X\omega\wedge \varphi$ is interpreted as improper integral $\lim_{\varepsilon\to 0}\int_{X\setminus V_{\varepsilon}}\omega\wedge\varphi$, where $V_{\varepsilon}$ is a tubular neighbhorhood of $Y$ of radius $\varepsilon$. The limit exists since $\frac{dz\wedge d\bar{z}}{z}=-2i\frac{dr\wedge d\theta}{e^{i\theta}}$ is integrable in a neighborhood of $z=0$. The residue map $\Res$ on the level of log-forms is defined as $\varphi\mapsto \mathcal{R}es_Y(\varphi)$, $\eta\wedge\frac{dz}{z}+\eta'\mapsto \eta$. We will study this in more details in Appendix \ref{Section_App}.

\begin{notation}\label{Notation_Pairing}
We denote the three pairings to be $\langle,\rangle_Y$, $\langle,\rangle^{X,Y}_U$ and $\langle,\rangle_X$ respectively.
\end{notation}

\begin{proposition}\label{MHSduality_prop} (i) The pairing on each column of the diagram \eqref{diagram} is unimodular (modulo torsion).\\
(ii) The diagram \eqref{diagram} is sign commutative. In fact, $\langle f(a),b\rangle_U^{X,Y}=-\langle a,q(b)\rangle_Y$ for all $a\in H_0^{2n-2}(Y,\Z)$ and $b\in H^{2n-1}(U,\Z)$.  $\langle c,p(d)\rangle_U^{X,Y}=\langle g(c),d\rangle_Y$ for all $c\in H^{2n-1}(X,Y,\Z)$ and $d\in H_0^{2n-1}(X,\Z)$.\\
(iii) The pairing induces isomorphisms between the dual mixed Hodge structures of the sequence \eqref{ext} twisted by $\mathbb Z(-2n+1)$ and the mixed Hodge structures on \eqref{ext_U}.
\end{proposition}
\begin{proof}
The argument is standard. The claim (i) dues to the natural pairing between the long exact sequence \eqref{les_XY} associated with the pair $(X,Y)$ and \eqref{les_U} to the pair $(X,U)$ in $\Z$-coefficent is unimodular. The unimodularity is preserved when taking associated kernels and cokernels.

(ii) can be verified from \eqref{RelativeOpenPair_Formula} directly. For (iii), it suffices to show that the dual mixed Hodge structure $H^{2n-1}(X,Y)^{\vee}$ is isomorphic to $H^{2n-1}(U)$. By definition of dual weight  \eqref{DualWeight_def} and Hodge filtrations \eqref{DualHodge_def}, it is equivalent to that the pairing induces an isomorphism
\begin{equation}
    W_{2n-1}H^{2n-1}(U)\cong (W_{2n-2}H^{2n-1}(X,Y))^{\vee} \label{weight_dual}
\end{equation}
on weight filtrations, and isomorphisms
\begin{equation}
    F^pH^{2n-1}(U,\C)\cong (F^{2n-p}H^{2n-1}(X,Y,\C))^{\vee} \label{Hodge_dual}
\end{equation}
on Hodge filtrations for $0\le p\le 2n-1$.

The isomorphism \eqref{weight_dual} follows from the commutativity of the pairing diagram \eqref{diagram}.
To show the isomorphism \eqref{Hodge_dual}, one concludes that the pairing
$$F^pH^{2n-1}(X,Y,\C)\times F^qH^{2n-1}(U,\C)\to \C$$
is zero for $p+q>2n-1$ by noticing the pairing  \eqref{RelativeOpenPair_Formula} vanishes on the level of cochains by type reason.
\end{proof}

\section{Curve Case}\label{Section_curve}

In this section, we study the relation between Abel-Jacobi map and the extension of mixed Hodge structures for dimension one case.

Let $D=\sum_i p_i-q_i$ be a degree zero divisor on a smooth projective curve $X$. Let's assume $p_i,q_j$ are disjoint for simplicity. The Abel-Jacobi map can be extracted from the information of mixed Hodge structure on $H^1(X\setminus|D|)$, where $|D|=\bigcup_i\{p_i,q_i\}$ is the support of the divisor. Indeed, there is an exact sequence of mixed Hodge structures
\begin{equation}
    0\to H^1(X,\Z)\to H^1(U,\Z)\xrightarrow{\textup{Res}} \Hv^0(|D|,\Z)(-1)\to 0,\label{curve_AJ_seq}
\end{equation}
where $\textup{Res}$ is the residue map, and $\Hv^0(|D|,\Z)$ is the $\Z$-submodule of $H^0(|D|,\Z)$ generated by $p_i-q_j$, $p_i-p_j$ and $q_i-q_j$ for all $i\neq j$. 

$D=\sum_ip_i-q_i\in \Hv^0(|D|)$ and defines a sub-exact sequence 
\begin{equation}
     0\to H^1(X,\Z)\to E_D\to\Z(-1)\to 0\label{Eqn_subseqCurve}
\end{equation}
of mixed Hodge structures of \eqref{curve_AJ_seq}. In particular, $E_D$ is a sub-mixed Hodge structure of $H^1(U)$.

The extension class $[E_D]$ of the exact sequence \eqref{Eqn_subseqCurve} can be obtained as follows. First, by the Riemann-Roch theorem, there is a meromorphic 1-form $\xi$ on $X$ whose poles are of the first order and are on $D$ and such that $\textup{Res}_{p_i}(\xi)=1$ and $\textup{Res}_{q_i}(\xi)=-1$. Note any other meromorphic 1-form $\xi'$ satisfying the same condition will differ by a holomorphic 1-form from $\xi$. Choose a basis $\{\gamma_1,\ldots,\gamma_{2g}\}$ for $H_1(X,\Z)$ and choose the presenting cycles disjoint from poles of $\xi$. Let $\eta\in C^{\infty}(X)$ be a smooth 1-form on $X$ such that $\int_{\gamma_i}\xi=\int_{\gamma_i}\eta$ for all $i=1,\ldots, 2g$ . By choosing a different integral basis $\{\gamma_1',\ldots, \gamma_{2g}'\}$ disjoint from the $p_i,q_i$, the resulting $\eta'$ will differ from $\eta$ by an element in $H^1(X,\Z)$. So $\eta$ is well-defined in
$$\textup{Ext}^1_{\textup{MHS}}(\Z(-1),H^1(X,\Z))\cong \frac{H^1(X,\C)}{F^1H^1(X,\C)+H^1(X,\Z)}$$
and coincides with the extension class $[E_D]$ of \eqref{Eqn_subseqCurve}.

Proposition \ref{Intro_Prop_Curve} is equivalent to the following statement.
\begin{proposition}\label{Prop_CurveDoubleInt}
For all $\omega \in H^0(X,\Omega_X)$,
 \begin{equation}
    \sum_i\int_{q_i}^{p_i}\omega=\int_X\omega\wedge \eta\mod \textup{periods}.\label{Intro_Eqn_iint}
 \end{equation}
\end{proposition}
\begin{proof}
First, there is a closed smooth 1-form $\varphi\in A^1_X(\log |D|)$ meromorphic around $p_i, q_i$ and satisfies that  (i) $\textup{Res}_{p_i}(\varphi)=1$, $\textup{Res}_{q_i}(\varphi)=-1$, (ii) $\varphi$ is supported on a simply connected domain containing $|D|$ (cf. \cite[Lemma 20.5]{Forster}). In particular on can choose $\gamma_i$ disjoint from $\varphi$, and choose $\eta=\xi-\varphi$. Now we can express $\int_X\omega\wedge \eta=\int_{X}\omega\wedge\xi-\int_X \omega\wedge\varphi$ and each term is interpreted as improper integral and exists since both $\xi$ and $\varphi$ have poles of order one. Now, $\omega\wedge\xi=0$ since $\xi$ is meromorphic and $\omega$ is holomorphic. Finally, $\int_{X}\varphi\wedge\omega=\sum_i\int_{q_i}^{p_i}\omega$ is due to the exactness of $\omega$ on the support of $\varphi$ and Stokes theorem \cite[Lemma 20.3]{Forster}.
\end{proof}

In the next section (Corollary \ref{Cor_Hodge}), we will prove a general argument using a similar proof idea. 

As a corollary, there is a Hodge theoretic version of Abel's theorem:
\begin{corollary}
The Abel-Jacobi image $A(D)$ is zero if and only if the extension \eqref{Eqn_subseqCurve} of the mixed Hodge structures splits over $\Z$.
\end{corollary}

\section{Zhao's Topological Abel-Jacobi and Real Splitting} \label{Section_ZhaoTAJ}

In this section, we review the definition of the topological Abel-Jacobi map and some basic properties in \cite{Zhao}. Further, we find the formula is actually a topological pairing between the relative cohomology and the relative homology. Using Lefschetz duality, we obtain a new formula \eqref{Eqn_ZhaoPairingResidue} for the topological Abel-Jacobi map and a simpler form for Hodge classes \eqref{Eqn_ZhaoPairingHodge}, which provides a proof of the main theorem in a special case. In the end, we relate Zhao's topological Abel-Jacobi to Deligne's $\R$-splitting of mixed Hodge structure on $H^{2n-1}(X,Y)$.

\subsection{Definition and Construction}

\begin{definition}\normalfont 
Zhao's topological Abel-Jacobi map \cite{Zhao} is a group homomorphism
$$A:\Hv^{2n-2}(Y,\Z)\to J_{\prim}(X)$$
that sends $\alpha\in \Hv^{2n-2}(Y,\Z)$ to the linear functional whose value on $[\omega]\in F^nH_{\prim}^{2n-1}(X)$ is
\begin{equation}
    A_{\alpha} ([\omega])=\int_{\Gamma}\omega-\int_{\gamma}\tau+\int_Yh_{\alpha}\wedge \tau+T(\omega),\mod \textup{periods.} \label{TAJ_LongFormula}
\end{equation}
\end{definition}

Here $\Gamma$ is a smooth $(2n-1)$-chain on $X$ whose boundary $\gamma$ is a smooth topological cycle that represents the Poincar\'e dual of the cohomology class $\alpha$. $\tau$ is a $(2n-2)$-form on $Y$ such that $d\tau=\omega_{|Y}$. Next, $h_{\alpha}$ is the harmonic representative of $\alpha$ with respect to the K\"ahler metric on $Y$ induced from $X$ (One can choose the K\"ahler metric on $X$ to be restriction of Fubini Study metric on $\mathbb P^N$).  Finally, $T$ is a differential current on $X$ such that $dT=-\int_{Y}h_{\alpha}\wedge(\cdot)$. So $T$ vanishes on $d^*$-closed forms.

\begin{proposition}
The map \eqref{TAJ_LongFormula} is well-defined, namely, it only depends on the cohomology classes of $\alpha$ and $[\omega]$.
\end{proposition}{}
\begin{proof}
First note that choosing a different smooth $(2n-1)$-chain $\Gamma'$ with $\partial\Gamma'=\gamma$ defines a linear functional differs by integration over $\Gamma-\Gamma'\in H_{2n-1}(X,\mathbb Z)_{\prim}$. 

Next, let $\omega'$ be another primitive form on $X$ such that $\omega-\omega'=du$, and $\tau'$ on $Y$ such that $d\tau'=\omega'_{|Y}$. Then $\tau=\tau'+u_{|Y}+dv$ for some $2n-3$ form $v$ on $Y$. Let $A_i$ denote the difference between the $i$-th term of \eqref{TAJ_LongFormula} applied to $(\tau,\omega)$ and $(\tau',\omega')$.  Then
\begin{equation*}
    A_1+A_2=\int_{\Gamma}du-\int_{\gamma}(u+dv)=\int_{\gamma}u-\int_{\gamma}u-\int_{\partial \gamma}v=0
\end{equation*}
by Stokes theorem and that $\gamma$ is closed. We also have 
\begin{equation*}
    A_3+A_4=\int_Yh_{\alpha}\wedge (u+dv)+T(du)=\int_Yh_{\alpha}\wedge u+dT(u)=0
\end{equation*}
by Stokes theorem and the fact that harmonic forms are closed and the definition of the current $T$. This shows that \eqref{TAJ_LongFormula} vanishes on the exact relative forms. 

\end{proof}

\subsection{Topological Abel-Jacobi Map is a Topological Pairing}

 To reduce the expression \eqref{TAJ_LongFormula}, we take $\omega_{h}$ as the harmonic representative of the class $[\omega]$. 
 
 \begin{proposition}\label{Prop_ddbarHarmonic}
Let $\omega_h$ be a harmonic form representing a class in $\Hp^{2n-1}(X,\C)$, then there exists a form $\sigma_h$ on $Y$ such that 
$$(\omega_{h})_{|Y}=dd^c\sigma_h.$$ 
\end{proposition}
\begin{proof}
Decompose $\omega_h=\sum_{p+q=2n-1}\omega_h^{p,q}$ into pure types. Then each $\omega_h^{p,q}$ is harmonic and therefore closed. Proposition \ref{Prop_ddbarXY} implies that there exists $\sigma_h^{p,q}$ such that $(\omega_{h}^{p,q})_{|Y}=dd^c\sigma_h^{p,q}$. We may take $\sigma_h=\sum\sigma^{p,q}_h$.
\end{proof}

\begin{proposition}
Zhao's topological Abel-Jacobi map $A_{\alpha}([\omega])$ modulo periods can be expressed as a topological pairing 
\begin{equation}
   H^{2n-1}(X,Y)\times H_{2n-1}(X,Y)\to \C,\label{Coh_Hom_XYpairing} 
\end{equation}
 \begin{equation}\label{TAJ2}
 \big((\omega_h,d^c\sigma_h), (\Gamma,\gamma)\big )\mapsto \int_{\Gamma}\omega_{h}-\int_{\gamma}d^c \sigma_h.
\end{equation}
\end{proposition} 
\begin{proof}
First, we take $\omega_{h}$ the harmonic representative of the class $[\omega]$. Then the fourth term of \eqref{TAJ_LongFormula} vanishes because harmonic forms are $d^*$-closed.
By Proposition \ref{Prop_ddbarHarmonic}, we may take $\tau=d^c \sigma_h$. Then the third term also vanishes because harmonic forms are $L^2$-orthogonal to $\textup{Im}(\partial)\oplus \textup{Im}(\bar{\partial})$. Therefore, Zhao's topological Abel-Jacobi map becomes the form \eqref{TAJ2}. This formula is obtained in \cite[Definition 2.1.2]{Zhao}. 

To show that \eqref{TAJ2} is topological, note that if $d(u,v)=(du,u_{|Y}-dv)$ is exact, then $\int_{\Gamma}du-\int_\gamma (u_{|Y}-dv)=\int_{\gamma}u-\int_{\gamma}u_{|Y}=0$. On the other hand, if $(\partial \Sigma+K,\partial K)$ is a relative boundary, where $\Sigma$ is a $2n$-chain in $X$ and $K$ is a $(2n-1)$-chain in $Y$, and $(u,v)$ is closed, then $\int_{\partial\Sigma+K}u-\int_{\partial K}v=\int_{K}u_{|Y}-\int_{K}dv=0$.

\end{proof}
As a consequence, the topological Abel-Jacobi map is independent of the choice of the K\"ahler metric of $X$. Also, the property \textbf{(P1)} in the introduction can be proved.
\begin{proposition} \cite[Proposition 2.1.1]{Zhao}
When $\alpha$ is an algebraic cycle, the topological Abel-Jacobi map \eqref{TAJ2} agrees with Griffiths's Abel-Jacobi map.
\end{proposition}
\begin{proof}
Note that $\int_{\gamma}d^c \sigma_h=\int_{\gamma}(d^c \sigma_h-id\sigma_h)=-2i\int_{\gamma}\partial\sigma_h$. By construction $\partial \sigma_h$ has type $(n,n-2)+(n+1,n-3)+\cdots$, whose integration along an analytic subvariety vanishes. Therefore in $(\ref{TAJ2})$, $A_{\alpha}([\omega])=\int_{\Gamma}\omega_h$, which coincides with Griffiths' definition \eqref{Griffiths'AJ}.
\end{proof}


\subsection{New Formula Using Lefschetz Duality}

We aim to obtain a different formula for the topological Abel-Jacobi map using the explicit pairing \eqref{RelativeOpenPair_Formula}  between $H^{2n-1}(X,Y)$ and $H^{2n-1}(U)$. Essentially this is to apply Lefschetz duality \eqref{Eqn_LD} to the formula \eqref{TAJ2}. 

We first need a lemma that lifts residue on cohomology to the level of differential forms with log poles. Let $\mathcal{R}es$ denote the residue on the level of differential forms with log poles on Y. Let $\textup{Res}$ denote the residue map on the level of cohomology.


\begin{lemma}\label{App_Lemma_Zlift}
    Suppose $\psi$ is a closed form on $Y$ representing a class $\alpha\in \Hv^{k-1}(Y,\Z)$. Let $\beta\in H^{k}(U,\C)$ be any class such that $\textup{Res}(\beta)=\alpha$, then there is a closed form $\varphi\in A^{k}_X(\log Y)$ such that $[\varphi_{|U}]=\beta$ and $\mathcal{R}es_Y(\varphi)=\psi$. 
\end{lemma}
\begin{proof}
There is an exact sequence of complexes
\begin{equation}
    0\to A_X'^*\to A_X^*(\log Y)\xrightarrow{\mathcal{R}es} A_Y^{*-1}\to 0. \label{Eqn_AcyclicComplex}
\end{equation}

$A_X'^*$ is the kernel of residue map $\mathcal{R}es$. The complex $A_X'^*$ contains $A_X^*$ as a subcomplex. The inclusion map $A_X^*\hookrightarrow A_X'^*$ is a quasi-isomorphism. The boundary homomorphism $\delta$ induces the Gysin homomorphism (\cite[p.212]{Voisin1}). It turns out that the long exact sequence of cohomology associated with \eqref{Eqn_AcyclicComplex} is $H^{k}(X)\to H^k(U)\to H^{k-1}(Y)\to H^{k+1}(X)$ whose truncation is the short exact sequence \eqref{ext_U} that we considered. First we choose any $\varphi'\in A_X^k(\log Y)$ representing the class $\beta$. Then $\mathcal{R}es_Y(\varphi')=\psi+d\psi_0$ for some $\psi_0\in A_Y^{k-2}$. Then choose any $\varphi_0\in A^{k-1}_X(\log Y)$ with $\mathcal{R}es_Y(\varphi_0)=\psi_0$, so $\varphi=\varphi'-d\varphi_0$ will satisfy the desired condition.
\end{proof}
Now, we take $k=2n-1$ in Lemma \ref{App_Lemma_Zlift}. Namely, for an integral vanishing cycle $\alpha\in \Hv^{2n-2}(Y,\Z)$, choose a closed form $\psi$ such that $[\psi]=\alpha$. Then choose any $\tilde{\alpha}_{\Z}\in H^{2n-1}(U,\Z)$ whose residue is $\alpha$. Then by Lemma \ref{App_Lemma_Zlift}, there is a closed log form $\tilde{\psi}\in A^{2n-1}_X(\log Y)$ that represents  $\tilde{\alpha}_{\Z}$ and such that $\mathcal{R}es_Y(\tilde{\psi})=\psi$. 
 Let $[\omega]\in F^n\Hp^{2n-1}(X)$ be of pure type $(p,q)$ and $\sigma$ being a $(p-1,q-1)$-form on $Y$ such that $\omega_{|Y}=dd^c\sigma$ (cf. Propositon \ref{Prop_ddbarXY}).
\begin{theorem}\label{Prop_NewPairing}
 Zhao's topological Abel-Jacobi map $A_{\alpha}([\omega])$ modulo periods can be expressed as a topological pairing
$$H^{2n-1}(X,Y)\times H^{2n-1}(U)\to \C,$$
\begin{equation}
 \big ( (\omega,d^c\sigma),-\tilde{\alpha}_{\Z}\big ) \mapsto -\int_X\omega\wedge \tilde{\psi}+\int_Yd^c\sigma\wedge \psi. \label{Eqn_ZhaoPairingResidue}
\end{equation}


\end{theorem}

\begin{proof}
By Poincar\'{e} duality and Lefschetz duality \eqref{Eqn_LD}, the exact sequence \eqref{ext_U} is identified with the exact sequence of homology of the pair $(X,Y)$. In other words, there is a diagram 
\begin{figure}[ht]
    \centering
\begin{tikzcd}
H^{2n-1}(U,\Z)\MySymb{dr}\arrow[d,"\cong","L.D."'] \arrow[r,"\Res"] &  H^{2n-2}(Y,\Z) \arrow[d,"\cong","P.D."']\\
H_{2n-1}(X,Y,\Z)\arrow[r,"\partial"] & H_{2n-2}(Y,\Z)
\end{tikzcd}
\end{figure}{}

\noindent which is commutative up to a sign, meaning that if $\beta\in H^{2n-1}(U,\Z)$, then  $P.D.(\textup{Res}(\beta))=-\partial L.D.(\beta)$. (Note this diagram factors through the first square of \cite[VI, Theorem 9.2]{Bredon}.)

In particular, the Lefschetz dual of the class $(\Gamma,\gamma)\in H_{2n-1}(X,Y,\Z)$ is a class $-\tilde{\alpha}_{\Z}\in H^{2n-1}(U,\Z)$ such that $\Res(\tilde{\alpha}_{\Z})=\alpha=P.D.^{-1}(\gamma)$. So the pairing $\langle (\omega,d^c\sigma), (\Gamma,\gamma)\rangle_{X,Y}$ in \eqref{Coh_Hom_XYpairing} coincides with the pairing $\langle (\omega,d^c\sigma), -\tilde{\alpha}_{\Z}\rangle_U^{X,Y}$ in \eqref{Eqn_ZhaoPairingResidue}.

Lemma \ref{App_Lemma_Zlift} guarantees the existence of a closed form $\tilde{\psi}'\in A^{2n-1}_X(\log Y)$ representing the class $\tilde{\alpha}_{\Z}$ such that $\mathcal{R}es(\tilde{\psi}')=\psi$. The formula \eqref{RelativeOpenPair_Formula} implies the pairing above equals $-\int_X\omega\wedge \tilde{\psi}'+\int_Yd^c\sigma\wedge \psi$. Note that any $\tilde{\psi}$ forestated differs from $\tilde{\psi}'$ by a closed form on $X$ whose class lies in $\Hp^{2n-1}(X,\Z)$, so the pairings induced by the two forms differ by $\int_{X}\omega\wedge(\tilde{\psi}-\tilde{\psi}')$ which belongs to the periods.
\end{proof}

\begin{remark}\normalfont
In the formula \eqref{Eqn_ZhaoPairingResidue} (also in \eqref{TAJ2}) for topological Abel-Jacobi map, the two terms together only depend on the (co)homology classes. Neither of each term separately is topological. For instance, note that $d^c\sigma$ is not even a closed form.
\end{remark}

For vanishing cycles that are Hodge classes, the topological Abel-Jacobi map will have a simpler form as below, which generalizes Proposition \ref{Prop_CurveDoubleInt} in curve case.

\begin{corollary}\label{Cor_Hodge} With the same notations as above, assume moreover that $\alpha\in \Hv^{2n-2}(Y,\Z)$ is a Hodge class represented by a closed form $\psi$ of type $(n-1,n-1)$. Then Zhao's topological Abel-Jacobi map on $\alpha$ modulo periods can be expressed as 
\begin{equation}
    -\int_{X}\omega\wedge \tilde{\psi}, \label{Eqn_ZhaoPairingHodge}
\end{equation}
where $\tilde{\psi}\in A^{2n-1}_X(\log Y)$ defines an integral class in $H^{2n-1}(U,\Z)$ and such that $\mathcal{R}es_Y(\tilde{\psi})=\psi$.
\end{corollary}
\begin{proof}
It suffices to show the second term of \eqref{Eqn_ZhaoPairingResidue} vanishes. Note $d^c\sigma=-2i\partial\sigma+id\sigma$. So Stokes theorem implies $\int_Yd^c\sigma\wedge\psi=-2i\int_Y\partial\sigma\wedge\psi$, which vanishes by type reason since $\omega$ has type $(p,q)$ with $p\ge n$ and $\partial\sigma$ has type $(p,q-1)$. 
\end{proof}

\begin{corollary}\label{Cor_Hodge2}
  Zhao's topological Abel-Jacobi map agrees with Schnell's topological Abel-Jacobi map on Hodge classes.
\end{corollary}
\begin{proof}
Since $\alpha\in \Hv^{2n-2}(Y,\Z)\cap F^{n-1}\Hv^{2n-2}(Y)$, there is a Hodge lift $\tilde{\alpha}_F\in F^nH^{2n-1}(U,\C)$ represented by a form $\varphi_F\in F^nA^{2n-1}_X(\log Y)$, so $\varphi_F-\varphi$ defines the extension class of \eqref{Eqn_HodgeClassExt}. Since the integral $\int_{X}\omega\wedge \varphi_F$ is well-defined and equals $0$ by type reason, $-\int_{X}\omega\wedge \varphi=\int_X\omega\wedge(\varphi_F-\varphi)$. So Zhao's topological Abel-Jacobi map coincides with Schnell's topological Abel-Jacobi map.
\end{proof}

\begin{remark}\normalfont

If $X$ is smooth projective of dimension $2n-1$ and an algebraic cycle $Z\in \mathcal{Z}^n(X)_{\textup{hom}}$ is (up to rational equivalence) contained in a smooth very ample hypersurface $Y$ of $X$, then Corollary \ref{Cor_Hodge2} will imply that Griffiths Abel-Jacobi map and Carlson's Abel-Jacobi map agree on primitive part of intermediate Jacobian (see Proposition \ref{Intro_Prop_CarGri}). For example, the condition forestated holds for all 1-cycles in a threefold: According to \cite[Theorem 5.8]{KleimanSmoothCycle} a 1-cycle $Z$ is rational equivalent to $Z'-mL$, where $Z'$ is smooth, and $L$ is a section by a general codimension two linear subspace and $m\in \Z$. In particular, $Z'$ and $L$ are disjoint. Strong Bertini theorem \cite{StrongBertini} guarantees the existence of $Y$. 
\end{remark}

\subsection{$\partial\bar{\partial}$-Lemma and Deligne's $\R$-Splitting of Mixed Hodge Structures}
For a class $[\omega]\in \Hp^{p,2n-1-p}(X)$, let $\sigma$ be a $(p-1,2n-2-p)$ form on $Y$ such that $\omega_{|Y}=dd^c\sigma$ obtained from Proposition \ref{Prop_ddbarXY}.

\begin{lemma}\label{RsectionXY_lemma}
The map
\begin{equation}\label{Eqn_RsectionXY}
s: \Hp^{2n-1}(X,\C)\to H^{2n-1}(X,Y,\C)
\end{equation}
that sends each $[\omega]\in \Hp^{p,2n-1-p}(X)$ to $ (\omega,d^c\sigma)$ defines a $\R$-linear section of \eqref{ext} and coincides with Deligne's $\R$-linear section.
\end{lemma}

\begin{proof} First of all, we show that the map is well-defined. Since $\omega_{|Y}=dd^c\sigma$, $(\omega,d^c\sigma)$ is closed in the mapping cone complex. If there is another $\sigma'$ such that $dd^c\sigma'=\omega_{|Y}$, then $d^c\sigma-d^c\sigma'$ is a $d$-closed form. However, it also lies in $\textup{Im}(\partial)\oplus \textup{Im}\bar{\partial}$, so it is $L^2$-orthogonal to the space of harmonic forms, so the cohomology class is zero. Therefore $s$ is a $\C$-linear map, and a section by definition.

Secondly, since $d^c$ is a real operator, $s$ descends to a $\R$-linear section $s_{\R}:\Hp^{2n-1}(X,\R)\to H^{2n-1}(X,Y,\R)$.

Finally, to show $s_{\R}$ is the Deligne's $\R$-linear section, by Proposition \ref{Prop_Rsection}, we need to show that $s$ preserves the Hodge filtration $F^p$. For $[\omega]$ of type $(p,q)$, $s([\omega])=(\omega,d^c\sigma)$ for a $(p-1,q-1)$-form $\sigma$ on $Y$. By subtracting an exact form $(0,id\sigma)$, we see that $(\omega,d^c\sigma)\sim (\omega,-2i\partial\sigma)\in F^pA_X^{2n-1}\times F^pA_Y^{2n-2}$, which implies $s(\omega)\in F^pH^{2n-1}(X,Y,\C)$ (cf. Proposition \ref{Prop_App_Fp}).
\end{proof}

\begin{remark}
For a class $[\omega]\in \Hp^{2n-1}(X)$ that has possibly mixed type, one may choose a harmonic representative $\omega_h$ and a form $\sigma_h$ on $Y$ such that $(\omega_{h})_{|Y}=dd^c\sigma_h$ as a result of Proposition \ref{Prop_ddbarHarmonic}. Then map \eqref{Eqn_RsectionXY} sends $[\omega]$ to $(\omega_h,\sigma_h)$.

\end{remark}

\section{Proof of the Main Theorem} \label{Section_Proof}

This section is devoted to proving Theorem \ref{main theorem}. As in Theorem \ref{Prop_NewPairing}, we already found the explicit form of Zhao's topological Abel-Jacobi map, which proves the second part of the theorem. It is left to prove that two definitions of topological Abel-Jacobi maps agree. The proof relies on duality properties of $\R$-split mixed Hodge structures and Lemma \ref{RsectionXY_lemma}. Note that we don't need the explicit integral formula \eqref{Eqn_ZhaoPairingResidue} in this section.

First, one note that $d^c$ is a real operator, so by integrating against forms in $\Hp^{2n-1}(X,\R)$, Zhao's topological Abel-Jacobi map \eqref{Eqn_ZhaoPairingResidue} takes values in the real torus $$J_{\prim}(X,\R):=\Hp^{2n-1}(X,\R)^{\vee}/H_{2n-1}(X,\Z)_{\prim}.$$ 

Now we need to identify the torus $J_{\prim}(X,\R)$ with $J_0(X,\R)$  defined in \eqref{Intro_Eqn_J0}.

By Proposition \ref{MHSduality_prop}, the pairing on the third column of \eqref{diagram} is unimodular, so the pairing induces isomorphisms $H^{2n-1}_0(X,\Z)\cong \Hp^{2n-1}(X,\Z)^{\vee}\cong H_{2n-1}(X,\Z)_{\prim}$, where the last isomorphism is induced by Poincar\'e pairing. It induces isomorphisms between real tori
\begin{equation}
  J_0(X,\R)=\frac{H^{2n-1}_0(X,\R)}{H^{2n-1}_0(X,\Z)}\cong \frac{\Hp^{2n-1}(X,\R)^{\vee}}{\Hp^{2n-1}(X,\Z)^{\vee}}\cong \frac{\Hp^{2n-1}(X,\R)^{\vee}}{H_{2n-1}(X,\Z)_{\prim}}=J_{\prim}(X,\R).  \label{J0_Jprim}
\end{equation}
So they are equivalent characterizations of the primitive real intermediate Jacobian of $X$.

We write $\omega$ for $[\omega]$ for the rest of the paper. To prove that Schnell's topological Abel-Jacobi map \eqref{realTAJ} coincides with Zhao's topological Abel-Jacobi map, we need to show for given classes $\alpha\in \Hv^{2n-2}(Y,\Z)$ and $\omega\in \Hp^{2n-1}(X,\R)$ the following equality holds modulo periods
$$A_{\alpha}(\omega)=\int_{X} \omega\wedge \big (s_{\R}^U(\alpha)-s_{\Z}^U(\alpha)\big ),$$
where $s_{\R}^U$ and $s_{\Z}^U$ is the Deligne's $\R$-linear section and a $\Z$-linear section of the sequence \eqref{ext_U}.

By Theorem \ref{Prop_NewPairing} and Lemma \ref{RsectionXY_lemma}, using the notations in  \ref{Notation_Pairing}, the main theorem is equivalent to show

\begin{proposition} \label{IdentitySchnellPair_prop}
For $\alpha\in \Hv^{2n-2}(Y,\Z)$ and $\omega\in H^{2n-1}(X,\R)$,
\begin{equation*}
   \langle s_{\R}^{X,Y}(\omega),s_{\Z}^U(\alpha)\rangle^{X,Y}_U=\langle \omega,(s_{\Z}'^{U}(\alpha)-s_{\R}^U(\alpha)\rangle_X \mod \textup{periods},
\end{equation*}
where $s_{\R}^{X,Y}$ is the Deligne's $\R$-linear sections of sequence \eqref{ext}, and $s_{\Z}'^U$ is another $\Z$-linear section of \eqref{ext_U}.
\end{proposition}

According to Proposition \ref{MHSduality_prop}, the sequence \eqref{ext_U} is identified with 
\begin{equation}
    0\to H^{2n-1}_0(X,\Z)^{\vee}\to H^{2n-1}(X,Y,\Z)^{\vee}\xrightarrow{-f^*} \Hv^{2n-2}(Y,\Z)^{\vee}\to 0
\end{equation}
twisted by $\Z(-2n+1)$ as an exact sequence mixed Hodge structures, where the third map has an extra minus sign due to the sign commutativity of the first square of \eqref{diagram}. This will tell us the relationship between Deligne's real splitting of the two sequences.

    We will first prove a linear algebra lemma. Suppose $E$ is an $\R$-split mixed Hodge structure that arises from an extension sequence
\begin{equation}
    \begin{tikzcd}\label{seq}
    \centering
    0\arrow{r} & V\arrow{r}{f} & E\arrow{r}{g} & W\arrow{r}\arrow[bend left=20]{l}{s_{\R}} & 0, 
\end{tikzcd}
\end{equation}
where $V$ and $W$ are pure Hodge structures of weight $w-1$ and $w$ respectively. $f$ and $g$ are morphisms of mixed Hodge structures of weight $0$. Let $s_{\R}:W\to E$ be the canonical Deligne's $\R$-linear section. Let 
\begin{equation}
    \begin{tikzcd}\label{dualseq}
    \centering
    0\arrow{r} & W^{\vee}\arrow{r}{g^*} & E^{\vee}\arrow{r}{-f^*}\arrow[bend left=33]{l}{s_{\R}^*} & V^{\vee}\arrow{r} & 0 
\end{tikzcd}
\end{equation}
be the dual sequence inducing the dual mixed Hodge structures, with weight $W^{\vee}$ and $V^{\vee}$ being $-w$ and $-w+1$. Both $g^*$ and $-f^*$ are morphisms of mixed Hodge structure of weight $0$.

\begin{lemma}\label{LA_lemma}
    The dual sequence \eqref{dualseq} is $\R$-split with the corresponding decomposition
\begin{equation}
    E^{\vee}\cong g^*W^{\vee}\oplus \ker(s_{\R}^*).\label{split}
\end{equation}
\end{lemma}
\begin{proof}
    First, the sequence \eqref{dualseq} is $\R$-split by Proposition \ref{Rsplit_prop}. Next, since $\ker(s_{\R}^*)\subseteq E^{\vee}$ is defined over $\R$ and is isomorphic to $V^{\vee}$ via $f^*$, it defines a splitting \eqref{split} of $E^{\vee}$ over $\R$. To show that it coincides with Deligne's $\R$-splitting, it suffices to show that the isomorphism of real vector space
    \begin{equation}
        -f^*|_{\ker (s_{\R}^*)}: \ker(s_{\R}^*)\to V^{\vee} \label{isoHS_eqn}
    \end{equation}
induces an isomorphism of the Hodge structures of weight $w$ after tensoring with $\C$

From now on, all vector spaces are complex. By abusing the notation, we omit the subscript $\C$.

It suffices to show that for each $l\in F^pV^{\vee}$, its inverse $\tilde{l}$ under $(\ref{isoHS_eqn})$ lies in $F^{p}\ker(s_{\R}^*)$. This is equivalent to show that for each $e\in F^{-p}E$, $\tilde{l}(e)=0$.

Via the splitting $E\cong f(V)\oplus s_{\R}(W)$, $e=e_1+e_2$, with $e_1\in f(V)$ and $e_2\in s_{\R}(W)$. Since Deligne's splitting $s_{\R}$ preserves Hodge filtration, $e_2=s_{\R}\circ g(e)$ lies in $F^{-p}E$ as well, so does $e_1=e-e_2$. Let $v\in V$ be the element such that $f(v)=e_1$, then $v\in F^{-p}V$ by strictness of $f$. By definition of $\tilde{l}$, it vanishes on the subspace $s_{\R}(W)$, so $\tilde{l}(e)=\tilde{l}(e_1)=\tilde{l}(f(v))=-f^*\tilde{l}(v)=l(v)$=0. The last equality is due to the assumption on $l$ and the fact that $V$ is a pure Hodge structure of weight $w-1$.
\end{proof}

Denote $s_{\R}^{\vee}$ (resp. $s_{\Z}^{\vee}$ and $s_{\Z}'^{\vee}$) the Deligne's $\R$-linear (resp. $\Z$-linear) section(s) of \eqref{dualseq}.  Denote $\langle,\rangle_E$ the natural pairing $E\times E^{\vee}\to \Z$, and similarly for $\langle,\rangle_W$.

Now, the proof of Proposition \ref{IdentitySchnellPair_prop} reduces to show 

\begin{proposition}
For each $\omega\in W$ and $\alpha\in V^{\vee}$, there is an equality 
  \begin{equation}
      \langle s_{\R}(\omega), s_{\Z}^{\vee}(\alpha) \rangle_E=\langle \omega, s_{\Z}'^{\vee}(\alpha)-s_{\R}^{\vee}(\alpha)\rangle_W\mod \textup{periods}.\label{equality}
  \end{equation}
\end{proposition}
\begin{proof}
  We take an integral basis $v_1,\ldots,v_k,w_1,\ldots,w_n$ of $E$, where $f^{-1}(v_i)$'s and $g(w_j)$'s form integral basis for $V$ and $W$ respectively. 
  We can find a real basis $u_1,\ldots, u_n$ of $s(W)$, such that 
  $$u_i=\sum_{j=1}^ka_{ij}v_j+w_i$$
  for some $a_{ij}\in \R$ and for each $i=1,\ldots,n$.
  
  On the other hand, $w_1^{\vee},\ldots,w_n^{\vee}, v_1^{\vee},\ldots,v_k^{\vee}$ form a dual basis of $E^{\vee}$. Moreover, $(g^*)^{-1}(w_j^{\vee})$'s form a basis of $W^{\vee}$ and is the dual basis of $g(w_i)$'s. By Lemma \ref{LA_lemma}, the decomposition $(\ref{split})$ is the Deligne's real splitting for $E^{\vee}$, so we can find a real basis $\rho_1,\ldots,\rho_n$ of the real subspace $\ker(s_{\R}^*)$ of $E^{\vee}$, with 
  $$\rho_j=-\sum_{i=1}^na_{ij}w_i^{\vee}+v_j^{\vee},$$
for $j=1,\ldots,k$.
  
 Finally, to check the equality $(\ref{equality})$, it suffices to check it on a basis. We take $\omega=g(w_i)$ and $\alpha=-f^*(v_j^{\vee})$, then $s_{\R}(\omega)=u_i$ and $s_{\Z}^{\vee}(\alpha)=v_j^{\vee}+\sum_kn_kw_k^{\vee}$, with $n_k\in \Z$. Then the left-hand side of \eqref{equality} is
 $$\langle s_{\R}(\omega),s_{\Z}^{\vee}(\alpha)\rangle_E=\langle u_i,v_j^{\vee}+\sum_kn_kw_k^{\vee}\rangle_E=a_{ij}+\textup{periods}.$$
 
On the other hand,  $s_{\R}^{\vee}(f^*(v_j^{\vee}))=\rho_j$ and $s_{\Z}'^{\vee}(\alpha)=v_j^{\vee}+\sum_km_kw_k^{\vee}$, with $m_k\in \Z$. So the right-hand side of \eqref{equality} is 
 \begin{align*}
    \langle \omega, s_{\Z}'^{\vee}(\alpha)-s_{\R}^{\vee}(\alpha)\rangle_W&=\langle g(w_i),\sum_{i=1}^n(a_{ij}+m_i)(g^*)^{-1}(w_i^{\vee})\rangle_W\\
     &=a_{ij}+\textup{periods},
 \end{align*}
 since $\langle g(w_i),(g^*)^{-1}(w_j^{\vee})\rangle=\delta_{ij}$.
 So both sides of the equality \eqref{equality} match up.
\end{proof}

\appendix
\section{An explicit pairing}\label{Section_App}

In this appendix, we aim to find an explicit formula \eqref{RelativeOpenPair_Formula} for the pairing
\begin{equation*}
    H^{2n-1}(X,Y,\Z)\times H^{2n-1}(U,\Z)\to \Z.
\end{equation*}

This pairing is obtained from a natural isomorphism \cite[Corollary B.14]{MHS}
\begin{equation}
   \Psi: H^{2n-1}(X,Y,\Z)\cong H^{2n-1}_c(U,\Z). \label{compactsupport_coh}
\end{equation}
together with the Poincar\'{e} pairing
\begin{equation}
   H^{2n-1}_c(U,\Z)\times H^{2n-1}(U,\Z)\to \Z,~([\psi],[\varphi])\mapsto \int_{X}\psi\wedge \varphi.\label{App_Eqn_UcUPair}
\end{equation}

We're going to establish the isomorphism \eqref{compactsupport_coh} explicitly. Let $A^*(X,Y)$ be the subcomplex of $A^*(X)$ consisting of $C^{\infty}$ forms on $X$ that vanish on $Y$. Its cohomology group is denoted as $H^*(X,Y)_G$, which is another model for relative cohomology \cite[Chapter XII, Theorem 3.1]{God}.

Our first step is to establish an explicit isomorphism between the cohomology theory $H^*(X,Y)_G$ \cite{God} and the relative de Rham cohomology $H^*(X,Y)$ obtained from the mapping cone complex. For preparation, let's take a sequence of tubular neighborhoods $T''\subseteq T'\subseteq T$ of $Y$, such that the closure of the previous one is contained in the next. Let $\pi: T\to Y$ be the projection. Then we choose a bump function $h$ that is smooth, takes constant value $1$ in $T''$, and is zero on $X\setminus T'$.

\begin{lemma}\label{App_Lem1}
The inclusion map $\psi\mapsto (\psi, 0)$ defines an isomorphism 
$$\Theta: H^k(X,Y)_G\xrightarrow{\cong} H^k(X,Y).\label{God-to-BT}$$

The inverse is given by
 \begin{equation}
  \Psi: (\omega,\tau)\mapsto \omega-d(\pi^*\tau\cdot h). \label{BT-to-God}  
\end{equation}
\end{lemma}

\begin{proof}

First, we show $\Psi$ is well-defined. If $(\omega,\tau)$ is closed, then $\omega_{|Y}=d\tau$ and the $\Psi(\omega,\tau)$ is closed and vanishes on $Y$. Also, if $(\omega,\tau)=d(\omega',\tau')$ is an exact form,  one has $$\Psi d(\omega',\tau')=\Psi(d\omega',\omega'_{|Y}-d\tau')=d (\omega'-\pi^*\omega'_{|Y}\cdot h)+d(\pi^*d\tau'\cdot h).$$ 

One notes that $\omega'-\pi^*\omega'_{|Y}\cdot h$ vanish on $Y$, and $\pi^*d\tau'\cdot h=d(\pi^*\tau\cdot h)\pm \pi^*\tau\wedge dh$. Since $dh$ vanishes in $T''$ and $dd(\pi^*\tau\cdot h)=0$, one concludes that  $\Psi d(\omega',\tau')=d\lambda$ for some $\lambda$ supported on $X\setminus T''$. So $\Psi$ sends exact forms to exact forms. (However, one should note that $\Psi$ is not well-defined on non-closed forms in general.)

Next, let’s show $\Psi$ is inverse to the map \eqref{God-to-BT}. Obviously, $\Psi\circ\Theta=Id$. On the other side, 
$$\Theta\circ\Psi(\omega,\tau)=\big (\omega-d(\tau\cdot h),0\big ),$$
which differ to $(\omega,\tau)$ by an exact form $d(\tau\cdot h,0)$, so $\Theta\circ\Psi=Id$ on the level of cohomology. Therefore, the two relative cohomology theories are isomorphic.
\end{proof}

Now a class in $H^k(X,Y)_G$ is represented by a closed $k$-form $\omega'$ vanishing on $Y$. Since $Y$ is homotopy equivalent to a tubular neighborhood $T$ of $Y$,  $\omega'_{|T}=du$ for a smooth $(k-1)$-form $u$ on $T$. Note we can choose $u$ to vanish on $Y$. Let $h$ be a bump function defined similarly as before.
\begin{lemma}(\cite[Chapter XII, Theorem 3.1]{God})\label{App_Lem2}
    There is an isomorphism
    $$H_c^{k}(U)\xrightarrow{\cong}H^k(X,Y)_G $$
    induced by the inclusion map. The inverse map $\chi$ is given by $\omega'\mapsto \omega'-d(u\cdot h)$. 
\end{lemma}

\begin{proposition}\label{ExplicitPairing_Prop}
The pairing \eqref{RelativeOpenPair} 

\begin{equation*}
     H^{2n-1}(X,Y,\Z)\times H^{2n-1}(U,\Z)\to \Z
\end{equation*}
is explicitly given by 
\begin{equation}
    ([\omega, \tau], [\varphi])\mapsto \int_{X}\omega\wedge \varphi-\int_Y\tau\wedge \mathcal{R}es_Y(\varphi). \label{pairing}
\end{equation}

Here $\varphi$ is a smooth form with logarithmic poles along $Y$, and the first integral is interpreted as an improper integral. The residue map $\mathcal{R}es$ sends each log form $\psi\wedge \frac{dz}{z}$ to $\psi$, where $z$ is a local equation of $Y$. 
\end{proposition}

\begin{proof}

Via the isomorphism $H^k_c(U)\cong H^k(X,Y)_G\cong H^k(X,Y)$ in Lemma \ref{App_Lem1} and Lemma \ref{App_Lem2}, we'll transform the pairing \eqref{App_Eqn_UcUPair} to \eqref{pairing}. Suppose we have $[(\omega,\tau)]\in H^{2n-1}(X)$, and $[\varphi]\in H^{2n-1}(U)$, where $\varphi$ has at worst a log pole along $Y$. We have the pairing \eqref{App_Eqn_UcUPair} equals to 
\begin{align*}
   \int_X\chi\circ\Psi(\omega,\tau)\wedge \varphi &=\int_X\big (\omega-d(\lambda\cdot h)\big )\wedge \varphi\\
   &=\int_{X}\omega\wedge \varphi-\int_{T'}d(\lambda\cdot h)\wedge \varphi,
\end{align*}
where $\lambda=\pi^*\tau+u$.
Note that the first term is interpreted as an improper integral, which is finite since $\varphi$ has a pole of order one.

Since $d(\lambda\cdot h)=d\lambda=\omega$ on $T''$, the second term becomes
$$\int_{T'\setminus T''}d(\lambda\cdot h)\wedge \varphi+\int_{T''}\omega\wedge\varphi.$$

By Stokes' theorem, the first term equals to $-\int_{\partial T''}\lambda\wedge \varphi$. 

Now we can assume $T''$ is a disk bundle of radius $\varepsilon$. Let $\varepsilon$ goes to zero, $\int_{T''}\omega\wedge\varphi\to 0$. Thanks to the residue theorem, we obtain

$$\lim_{\varepsilon\to 0}\int_{\partial T''}\lambda\wedge \varphi=\int_{Y}\lambda\wedge \mathcal{R}es(\varphi),$$
using the fact that $u$ vanishes on $Y$, which establishes \eqref{pairing}.
\end{proof}

\section{Hodge filtration on relative cohomology}\label{Sec_AppB}
In this appendix, we will study the mixed Hodge structure for the relative cohomology $H^*(X,Y)$. We first review a sheaf theoretical definition of relative cohomology in \cite[p.78]{MHS} and apply it to a geometric case.

Let $i:Y\hookrightarrow X$ be the inclusion of compact K\"ahler manifolds. The $\Z$-coefficient relative cohomology is defined to be the hypercohomology of the mapping cone complex of the restriction of constant sheaf $i^*:\underline{\Z}_X\to \underline{\Z}_Y$. 
\begin{equation*}
H^k(X,Y,\Z):=\mathbb H^k(X,\underline{\mathcal{C}one}^{\bullet}(i^*:\underline{\Z}_X\to \underline{\Z}_Y)).
\end{equation*} 
The exact sequence of complexes of sheaves
$$0\to \underline{\Z}_Y\to \underline{\mathcal{C}one}^{\bullet}(i^*:\underline{\Z}_X\to \underline{\Z}_Y)) \to \underline{\Z}_X[1]\to 0$$
induces the long exact sequence of cohomology
\begin{equation}\label{Eqn_App_les}
\cdots \to H^{k-1}(Y,\Z)\to H^k(X,Y,\Z)\to H^k(X,\Z)\to \cdots,
\end{equation}
which is the long exact sequence of the cohomology of the pair $(X,Y)$.
 
$H^k(X,Y,\Z)$ carries a mixed Hodge structure that makes \eqref{Eqn_App_les} exact in the category of mixed Hodge structures. In fact, since the de Rham complex $\Omega_X^{\bullet}$ is a resolution of the constant sheaf $\underline{\C}_X$, which is the complexification of $\underline{\Z}_X$. The $\C$-coefficient relative cohomology is 
\begin{equation*}
H^k(X,Y,\C)=\mathbb H^k(X,\underline{\mathcal{C}one}^{\bullet}(i^*:\Omega_X^{\bullet}\to \Omega_Y^{\bullet})).
\end{equation*} 

The Hodge filtration $F^pH^k(X,Y,\C)$ is defined to be the hypercohomology of the complex of sheaves \cite[p.76, Theorem 3.22, and p.76]{MHS}
\begin{equation}\label{Eqn_OmegaXY}
\underline{\mathcal{C}one}^{\bullet}(i^*:F^p\Omega_X^{\bullet}\to F^p\Omega_Y^{\bullet}),
\end{equation}
where $F^p\Omega_X^{\bullet}=\Omega_X^{\ge p}$ is the $p$-th naive filtration of the holomorphic de Rham complex.

\begin{proposition}\label{Prop_App_Fp}

$F^pH^k(X,Y,\C)$ coincides with the $k$-th cohomology of the complex 
\begin{equation}\label{Eqn_App_ConeFpA}
\Cone^{\bullet}(i^*:F^pA_{X}^{\bullet}\to F^pA_{Y}^{\bullet}),
\end{equation}
where $A_X^{\bullet}$ is the complex of global $\C$-valued smooth forms on $X$ and $F^pA_{X}^{\bullet}=\oplus_{i\ge p,j}A_X^{i,j}$ is the subcomplex of sheaves consisting of smooth forms with at least $p$ holomorphic differential. 

\end{proposition}
\begin{proof}
The complex $\underline{\mathcal{C}one}^{\bullet}(i^*:\Omega_X^{\bullet}\to \Omega_Y^{\bullet})$ admits a fine resolution by a double complex
\begin{figure}[ht]
    \centering
    \begin{equation}\label{Eqn_App_diagram}
\begin{tikzcd}
&\cdots\arrow[r] &\mathcal{A}_X^{k,\bullet}\oplus\mathcal{A}_Y^{k-1,\bullet}\arrow[r,"d"] &  \mathcal{A}_X^{k+1,\bullet}\oplus\mathcal{A}_Y^{k,\bullet} \arrow[r]&\cdots\\
\underline{\mathcal{C}one}^{\bullet}(i^*\text{:}\Omega_X^{\bullet}\to \Omega_Y^{\bullet})\arrow[r,symbol=\text{=}]&\big\{\cdots\arrow[r]&\Omega_X^k\oplus \Omega_Y^{k-1}\arrow[r,"d"] \arrow[u,hook]& \Omega_X^{k+1}\oplus \Omega_Y^{k}\arrow[u,hook]\arrow[r]&\cdots\big\}.
\end{tikzcd}
\end{equation}
\end{figure}{}

\noindent with $\mathcal{A}_X^{k,j}$ is the sheaf of $\C$-valued smooth $(k,j)$ forms on $X$ and $\mathcal{A}_X^{k,\bullet}$ is the complex $\{\cdots\to\mathcal{A}_X^{k,j}\xrightarrow{\bar{\partial}}\mathcal{A}_X^{k,j+1}\to\cdots\}$. Then
the total complex of the double complex is just the mapping cone complex of $i^*: \mathcal{A}_X^{\bullet}\to \mathcal{A}_Y^{\bullet}$. Apply $F^p$ to the resolution diagram \eqref{Eqn_App_diagram}, one obtains a fine resolution of \eqref{Eqn_OmegaXY} and concludes that it is quasi-isomorphic to the mapping cone complex of $i^*: F^p\mathcal{A}_X^{\bullet}\to F^p\mathcal{A}_Y^{\bullet}$, whose hypercohomology is just cohomology of the complex of the global sections \eqref{Eqn_App_ConeFpA}.
\end{proof}

\bibliographystyle{plain}
\bibliography{bibfile}

\end{document}